\documentclass{amsart}

\newtheorem{theorem}{Theorem}[section]

\newtheorem{proposition}[theorem]{Proposition}
\newtheorem{corollary}[theorem]{Corollary}

\theoremstyle{definition}

\newtheorem{remark}[theorem]{Remark}
\newtheorem*{acknowledgement}{Acknowledgement}

\theoremstyle{remark}

\newcommand\mylabel[1]{\label{#1}}

\newcommand{\KK}{\mathbb{K}}

\newcommand{\ZZ}{\mathbb{Z}}

\newcommand  {\shB}     {\mathcal{B}}

\newcommand  {\shExt}   {\mathcal{E} \!\text{\textit{xt}}}
\newcommand  {\shE}     {\mathcal{E}}
\newcommand  {\shF}     {\mathcal{F}}

\newcommand  {\shH}     {\mathcal{H}}
\newcommand  {\shHom}   {\mathcal{H}\!\text{\textit{om}}}
\newcommand  {\shI}     {\mathcal{I}}

\newcommand  {\shM}     {\mathcal{M}}

\newcommand  {\shN}     {\mathcal{N}}
\newcommand  {\shL}     {\mathcal{L}}

\newcommand  {\shS}     {\mathcal{S}}
\newcommand  {\shT}     {\mathcal{T}}
\newcommand  {\shTor}   {\mathcal{T}\!\text{\textit{or}}\,}

\newcommand  {\foU}     {\mathfrak{U}}
\newcommand  {\foV}     {\mathfrak{V}}

\newcommand  {\foX}     {\mathfrak{X}}

\newcommand  {\Ass}     {\operatorname{Ass}}

\newcommand  {\depth}   {\operatorname{depth}}

\newcommand  {\Ext}     {\operatorname{Ext}}

\newcommand  {\GL}      {\operatorname{GL}}

\newcommand  {\Grass}   {\operatorname{Grass}}

\newcommand  {\Hom}     {\operatorname{Hom}}

\newcommand  {\lra}     {\longrightarrow}

\newcommand  {\maxid}   {\mathfrak{m}}

\newcommand  {\naive}   {{\operatorname{naive}}}
\newcommand  {\nor}     {{\operatorname{nor}}}

\newcommand  {\primid}  {\mathfrak{p}}
\renewcommand{\O}       {\mathcal{O}}

\newcommand  {\perf}    {{\operatorname{perf}}}
\newcommand  {\Pic}     {\operatorname{Pic}}
\newcommand  {\pd}      {\operatorname{pd}}

\newcommand  {\ra}      {\rightarrow}

\newcommand  {\red}     {{\operatorname{red}}}

\newcommand  {\Sing}    {\operatorname{Sing}}

\newcommand  {\Spec}    {\operatorname{Spec}}

\def\mydate{\number\day\space\ifcase\month \or January\or February\or March\or
April\or May\or June\or July\or
August\or September\or October\or November\or December\fi \space\number\year}

\usepackage{amscd,amssymb}
\usepackage[PostScript=dvips]{diagrams}

\begin{document}

\title[Vector bundles and global resolutions]
         {Existence of vector bundles and global resolutions
          for singular surfaces}

\author[Stefan Schroer]{Stefan Schr\"oer}
\address{Mathematische Fakult\"at, Ruhr-Universit\"at,
            44780 Bochum, Germany}
\email{s.schroeer@ruhr-uni-bochum.de}

\author[Gabriele Vezzosi]{Gabriele Vezzosi}
\address{Dipartimento di Matematica, Universit\`a di Bologna,
            40127 Bologna, Italy}
\email{vezzosi@dm.unibo.it}


\dedicatory{Final version, 26 April 2002}

\begin{abstract}
We prove two results about vector bundles
on  singular algebraic surfaces. First, on \emph{proper surfaces} there
are vector bundles of rank two with arbitrarily large second Chern number
and fixed determinant.
Second, on \emph{separated normal surfaces} any coherent sheaf
is the quotient of a vector bundle.
As a consequence, for such surfaces the Quillen $K$-theory of vector
bundles  coincides with the Waldhausen $K$-theory of perfect
complexes.
Examples show that, on nonseparated schemes,   usually many
coherent sheaves   are not quotients of vector bundles.
\end{abstract}

\maketitle

\tableofcontents

\section*{Introduction}

Let $Y$ be a  scheme. We pose the following question:
Does there exist a nontrivial locally free $\O_Y$-module of
finite rank?
Many fundamental problems in algebraic geometry are
related to this.
To mention a few:  in the theory of \textit{moduli spaces}, one
wants to know whether certain moduli spaces
of semistable vector bundles with given Chern numbers
are nonempty. In $K$-\textit{theory}, one seeks to replace coherent
sheaves
by complexes of locally free sheaves, or perfect complexes
by bounded complexes of locally free sheaves.
In the theory of \textit{Brauer groups}, a fundamental problem is to find
a locally free sheaf admitting an Azumaya
algebra structure with given cohomology class.

There is a good chance to tackle such problems if $Y$ admits
an ample invertible sheaf, or more generally an ample family of
invertible sheaves (\cite{SGA 6}, Expos\'e II).
Very few things, however, seem to be known in general. For example,
it is open whether or not a proper algebraic scheme has
nontrivial vector bundles at all.

The goal of this paper is to attack the problem for the
simplest family
of schemes where few results are known, namely algebraic surfaces
with
singularities over an arbitrary ground field.
Our first main result (Theorem \ref{existence}) is that a
\textit{proper} algebraic surface admits rank two vector bundles
with arbitrary preassigned determinant and arbitrarily large second
Chern number. The second main result (Theorem \ref{quotient}) is
that on \textit{normal separated} algebraic surfaces, any coherent
sheaf is the quotient of a locally free sheaf. In other words, any
coherent sheaf admits a global resolution with locally free sheaves.
Examples show that this fails for trivial reasons without
the separatedness assumption (Proposition \ref{vector bundles II}).
Using Theorem \ref{quotient}
we infer that on   normal separated   surfaces
Quillen's $K$-theory of vector bundles
coincides with the Waldhausen $K$-theory of perfect complexes
(Theorem \ref{K-spectra and higher K-groups}).

The situation for algebraic surfaces with singularities
is similar to the case of \emph{nonalgebraic} smooth compact complex
surfaces. Schuster \cite{Schuster 1982} extended   existence
of global resolution for smooth algebraic surface to smooth
compact complex surfaces. For nonalgebraic smooth compact
complex surfaces, B\u{a}nic\u{a} and Le Potier
\cite{Banica; Le Potier 1987} proved that    Chern numbers for
torsion free sheaves satisfy certain inequalities. For rank two
vector bundles, this reduces to $2c_2\geq c_1^2$, so the second
Chern number is bounded from below if the determinant is fixed.

This analogy  breaks down in higher dimensions.
Burt Totaro pointed out to us recent results of
Claire Voisin: among other things, she proved   that
ideal sheaves of closed points on general complex
tori of dimension $\geq 3$ are not quotients of
locally free sheaves (\cite{Voisin 2001}, Corollary 2).

\begin{acknowledgement}
We thank the referee for his remarks and for pointing
out some errors in the first version of the paper.
\end{acknowledgement}
\section{Existence of vector bundles}

Throughout this paper we work over a fixed (but arbitrary) ground
field $k$. In this section, $Y$ denotes a \textit{proper} surface.
The word \emph{surface} means a $k$-scheme of finite
type whose irreducible components are 2-dimensional.
We pose the following question: ``given an invertible
$\O_Y$-module $\shL$ and an integer $c_2\in\ZZ$, does there
exist a locally free $\O_Y$-module $\shF$ of rank two
with determinant $\det(\shF)\simeq\shL$ and second
Chern number $c_2(\shF)=c_2$?''.

Schwarzenberger proved such a result
for \textit{smooth} surfaces over algebraically closed ground fields
(\cite{Schwarzenberger 1961}, Theorem 8).
Beyond this, however, little seems to be known.
It might easily happen that $\Pic(Y)=0$ for normal surfaces
(see \cite{Schroeer 1999}, Section $3$), and it is a priori
unclear whether or not there are any nontrivial vector bundles on
arbitrary singular surfaces (though, by \cite{Schroeer 2001},
Theorem 4.1, any proper normal surface has nontrivial vector
bundles of square rank).

To start with, let us briefly recall the notions of $S_1$-ization and
$S_2$-ization since we will use them quite often.
Given (a surface) $Y$, its $S_1$-ization is the closed subscheme
$Y'\subset Y$ having the same underlying topological space as
$Y$ and whose defining ideal is the ideal $\shI\subset\O_Y$ of
sections whose support has dimension $\leq 1$.
The $S_2$-ization $f:X\ra Y$  of $Y$ is defined as follows:
For any $y\in Y$, the stalk $f_*(\O_X)_y$ is the intersection
of all local rings $\O_{Y',x}$ for all points $x\in \Spec(\O_{Y,y})$ of
codimension one, the intersection taking place in the ring of
total fractions for $\O_{Y',y}$. According to
\cite{EGA IVb}, Proposition $5.11.1$, the
$S_2$-ization $f:X\ra Y$ is a finite birational morphism.
Note that $Y'$ and $X$ satisfy Serre's condition $(S_1)$ and
$(S_2)$, respectively.

Next, it is perhaps a good idea to define  for our situation  what we
mean by the \emph{second Chern number} $c_2(\shF)\in\ZZ$ for a given
locally free $\O_Y$-module $\shF$ of rank $r\geq 0$. The following
naive approach suits us well. Let $f:X\ra Y$ be
the $S_2$-ization. To define $c_2(\shF)$, set $\shE=f^*(\shF)$ and
$\shL=\det(\shE)$. The Riemann--Roch equation
\begin{equation}
\label{RR}
\chi(\shE)=
\frac{1}{2} \shL\cdot(\shL-\omega_X) -c_2(\shE) + r\chi(\O_X)
\end{equation}
\emph{defines} an integer $c_2(\shE)$, and we set
$c_2(\shF)=c_2(\shE)$. Here $\omega_X$ is the dualizing sheaf, which
exists because $X$ is Cohen--Macaulay.

We now state our first main result:

\begin{theorem}
\mylabel{existence}
Let $Y$ be a proper surface. Then for each invertible
$\O_Y$-module $\shL$  and each integer $c_2\geq 0$, there is
a locally free
$\O_Y$-module
$\shF$ of rank two satisfying $\det(\shF)\simeq\shL$ and
$c_2(\shF)\geq c_2$.
\end{theorem}

The proof requires some preparation. Suppose we have a proper
birational morphism $f:X\ra Y$. The union of all integral
curves mapped to a point is called the \emph{exceptional curve}
$E\subset X$. The union of all integral curves
given by the codimension one points $x\in X$ for which
$\O_{Y,f(x)}\ra\O_{X,x}$ is finite but not bijective is called the
\emph{ramification curve}. The following tells us that in
constructing vector bundles, we may pass to certain auxiliary
surfaces.

\begin{proposition}
\mylabel{bijection}
Let $f:X\ra Y$ be a proper birational morphism with empty
ramification curve. Suppose the
exceptional curve $E\subset X$ supports an effective
Cartier divisor $D\subset X$ with $\O_E(-D)$ ample.
Then for each $r\geq 0$, there is an $n\geq 1$ such
that the functor $\shF\mapsto f^*(\shF)$ induces a surjection
from the set of isomorphism classes of
locally free $\O_Y$-modules of rank $r$ to
the set of isomorphism classes of locally free
$\O_X$-modules of rank $r$ whose restriction to $nD$ is trivial.
\end{proposition}

\proof
First consider the special case that $f:X\ra Y$ is finite.
Then the exceptional curve $E\subset X$ is empty, and
$\O_Y\ra f_*(\O_X)$ is bijective outside finitely many closed points.
Fix a locally free $\O_X$-module $\shE$ of rank $r$.
For each $y\in Y$, the preimage $f^{-1}(\Spec\O_{Y,y})$ is a
semilocal scheme, on which $\shE$ becomes trivial. Consequently,
there is an affine open covering $V_i\subset Y$ such that $\shE$
becomes trivial on the preimages $U_i=f^{-1}(V_i)$. Set $\foV=(V_i)$
and
$\foU=(U_i)$. Passing to a refinement of $\foV$, we may assume that
each
$V_i$ contains at most one point $y\in Y$ where $\O_Y\ra f_*(\O_X)$
is not bijective. Then $U_i\cap U_j=V_i\cap V_j$ for $i\neq j$, hence
the canonical map on alternating cocycles
$Z^1(\foV,\GL_{r,Y})\ra Z^1(\foU,\GL_{r,X})$ is bijective. It follows
that the functor
$\shF\mapsto f^*(\shF)$ is surjective on isomorphism classes.

Now let $f:X\ra Y$ be arbitrary.
Applying the preceding special case to the Stein factorization
$\Spec(f_*\O_X)\ra Y$, we may replace $Y$ by $\Spec(f_*\O_X)$ and
assume that
$\O_Y\ra f_*(\O_X)$ is bijective. By the Theorem on Formal Functions
(\cite{EGA IIIa}, Theorem 4.1.5), the functor $\shF\mapsto f^*(\shF)$
is an equivalence between the categories of locally free
$\O_Y$-modules $\shF$ of rank $r$ and the category of locally free
$\O_X$-modules
$\shE$ of rank $r$
whose restriction to the formal completion $\foX=X_{/E}$ is trivial.
The same argument as in \cite{Schroeer 2001}, Lemma 2.2, shows that
there is an integer $n\geq 1$
such that $\shE$ is trivial on $\foX$ if and only if it is trivial
on
$nD$, and the result follows.

Note that  \cite{Schroeer 2001}, Lemma 2.2 is stated for
resolutions of singularities of normal surfaces, but the proof holds
literally true in our situation. That result also involves
a family $\shB$ of locally free $\O_E$-modules that is bounded up to
tensoring with line bundles; here we apply the result with
the trivial family $\shB=\left\{\O_E\right\}$.
\qed

\begin{remark}
\mylabel{stein}
If $\O_Y\ra f_*(\O_X)$ is bijective, the preceding arguments show
that the functor $\shF\mapsto f^*(\shF)$
yields a bijection between the set
of isomorphism classes of locally free $\O_Y$-modules of rank $r$
and the set of isomorphism classes of locally free
$\O_X$-modules of rank $r$ whose restriction to $nD$ is trivial.
\end{remark}

Passing to such auxiliary surfaces $X$ does not change
second Chern numbers:

\begin{proposition}
\mylabel{number}
Let $f:X\ra Y$ be a proper birational morphism, $R\subset X$ its
ramification curve, $\shF$ a locally free $\O_Y$-module
of finite rank, and $\shE=f^*(\shF)$ its pullback. If
$\shF$ is trivial on $f(R)$, then $c_2(\shF)=c_2(\shE)$.
\end{proposition}

\proof
Making base change with the $S_1$-ization of $Y$ and replacing
$X$ by its $S_2$-ization, we may assume that $Y$ has no embedded
components and that $X$ is Cohen--Macaulay.
If $f:X\ra Y$ is finite, it must factor over the $S_2$-ization of
$Y$. If $\O_Y\ra f_*(\O_X)$ is bijective, then $Y$ already satisfy
Serre's condition
$(S_2)$. Using Stein factorization, we therefore may assume that $Y$
is Cohen--Macaulay.

The trace map $f_*(\omega_X)\ra\omega_Y$ is
bijective outside the union of $f(R)$ and an additional finite
subset. Since
$\shM=\det(\shF)$ is trivial on
this set, the projection formula gives $\shM\cdot\omega_Y=
f^*(\shM)\cdot \omega_X$.
Furthermore, we have
$\chi(\shF)-\chi(\shE)=r\chi(\O_Y)-r\chi(\O_X)$, because $\shF$ is
trivial on $f(R)$.
Now the Riemann--Roch formula (\ref{RR})
immediately implies
$c_2(\shF)=c_2(\shE)$.
\qed

\medskip
Now we come to the  construction of auxiliary surfaces:

\begin{proposition}
\mylabel{blow-up}
Let $Y$ be a proper surface. Then there is
a projective birational morphism $f:X\ra Y$ with empty
ramification curve such that the following holds:
\begin{enumerate}
\renewcommand{\labelenumi}{(\roman{enumi})}
\item There is a Cartier divisor $D\subset X$
supported by the exceptional curve $E\subset X$
with $\O_E(-D)$ ample;
\item There is an irreducible Cartier divisor $C\subset X$
contained in the Cohen--Macaulay locus of $X$
and having no  component in common with $E\subset X$.
\end{enumerate}
\end{proposition}

\proof
Choose a codimension one point $y\in Y$ not contained in
$\Ass(\O_Y)$ so that $Y_\red$ is normal near $y\in Y_\red$.
There is a regular
element
$s\in\maxid_{Y,y}$ such that
$\O_{Y,y}/(s)$ is 0-dimensional. Let $B\subset Y$ be
the schematic image of the canonical morphism $\Spec\O_{Y,y}/(s)\ra
Y$. Then
$B\subset Y$ is an irreducible curve and Cartier outside
finitely many closed points.

We now make a sequence of blow-ups $Y_4\ra\ldots\ra Y_1\ra Y$
whose centers are infinitesimal near to $B$ and have $0$-dimensional
image on $B$. We shall denote by $B_i\subset Y_i$ the strict
transform of
$B$. First, we blow up the intersections
$B\cap\overline{\left\{y\right\}}$, where $y$ ranges over
the nongeneric associated points   $y\in\Ass(\O_X)$.
Then the strict transform $B_1\subset Y_1$   becomes disjoint
from   embedded components. Next, let $H_1\subset Y_1$ be the branch
locus for the normalization map $Y_1^\nor\ra Y_1^\red$, and
$Y_2\ra Y_1$ be the blowing up with center
$H_1\cap B_1$.
Write $R_2\subset Y_2$ for the exceptional curve of $Y_2\ra Y_1$,
and $H_2$ for the strict transform of $H_1$.

To define the next blow up, we first choose an   open
subset $U_2\subset Y^\red_2$ containing the support of $B_2$
and disjoint
from $H_2$, so that $R_2\cup B_2$ contains the singular locus of
$U_2$. Let
$r:U_3\ra U_2$ be a resolution of singularities. Then there is an
effective Cartier divisor
$R_3\subset U_3$ supported  by $r^{-1}(\Sing(U_2))$ such that
$-R_3$ is relatively ample with respect to $r$. Consequently,
$r_*\O_{U_3}(-nR_3)\subset r_*\O_{U_3}$ is contained in
$\O_{U_2}$ for
$n\gg 0$, and it then follows as in the proof of
\cite{EGA IIIa}, Proposition 2.3.5  that $r:U_3\ra U_2$
is the blowing up of some center $Z_2\subset \Sing(U_2)$.
Now let $Y_3\ra Y_2$ be the blowing up of the schematic closure
$\overline{Z_2}\subset Y_2$. By construction, the strict transform
$B_4\subset Y_4$ is supported by the regular locus of $Y_4^\red$.

Let $g:Y'\ra Y$ be the composite morphism for the preceding sequence
of blowing-ups. According to \cite{Raynaud; Gruson 1971}, Lemma
5.1.4, this is a blowing-up of a 0-dimensional closed subscheme
$Z\subset Y$ supported by $B$. The preimage $D'=g^{-1}(Z)$ is an
effective Cartier divisor supported on the exceptional curve, whose
ideal is isomorphic to $\O_{Y'}(1)$. In other words, it is ample
on the exceptional curve.

Let $B'\subset Y'$ be the strict transform of $B$.
Let $h:X\ra Y'$ be
the blowing-up of
$B'$. This is a homeomorphism, because $B'\cap Y'_\red\subset
Y'_\red$ is already a Cartier divisor. Consequently, the effective
Cartier divisor
$C=h^{-1}(B')$ has no irreducible component in common with the
exceptional curve $E\subset X$ for the composition
$f:X\ra Y$. Replacing $X$ in some open neighborhood of
$C$ by its $S_2$-ization, we may also assume that $C$ lies in the
Cohen--Macaulay locus. Since $h$ is finite, the Cartier divisor
$D=h^{-1}(D')$ has $\O_E(-D)$ ample.
\qed

\medskip
\emph{Proof of Theorem \ref{existence}:}
First choose a proper birational morphism $f:X\ra Y$ and
Cartier divisors $C\subset X$ and $D\subset X$ as in Proposition
\ref{blow-up}. Then choose an integer $n\geq 1$
as in Proposition \ref{bijection} so that any locally free
$\O_X$-modules of rank two that is trivial on $nD$ comes from a
locally free
$\O_Y$-modules of rank two .

In light of Proposition \ref{number}, we have to solve the
following problem: Given an invertible $\O_X$-module $\shL$ that is
the preimage of an invertible $\O_Y$-module, find a locally free
$\O_X$-module $\shE$ of rank two whose restriction to $nD$ is
trivial, so that
$\det(\shE)\simeq \shL$ and
$c_2(\shE)\gg 0$.
The idea is to start with $\shE_0=\shL\oplus\O_X$ and apply two
elementary transformations (see \cite{Maruyama 1982})
that \textit{cancel each other on} $nD$.
Note that $\det(\shE_0)=\shL$ and $c_2(\shE_0)=0$.

To proceed, choose an ample effective Cartier divisor $A\subset C$
and an integer $t>0$ so that $H^1(C,\shL_C((t-1)A))=0$.
Let $s_1:\shL_C\ra\shL_C(tA)$ be the map induced by $tA$.
The exact sequence
$$
H^0(C,\shL_C(tA))\lra H^0(A,\shL_A(tA))\lra H^1(C, \shL_C((t-1)A))
$$
shows that some map $s_2:\O_C\ra\shL_C(tA)$ has no zeros
on $A$. This gives a surjection
$s:\shE_0|_C\ra \shL_C(tA)$. In turn, the exact sequence
\begin{equation}
\label{sequence}
0\lra\shE_1\lra \shE_0\lra \shL_C(tA)\lra 0
\end{equation}
defines a coherent $\O_X$-module $\shE_1$ of rank two, which is
locally free because the stalks of $\shL_C(tA)$ have projective
dimension one. In other words, $\shE_1$ is an elementary
transformation of $\shE_0$. Taking determinants on
$X-A$, we compute
$\det(\shE_1)=\shL(-C)$. This works because $X$ satisfies Serre's
condition
$(S_2)$ near $A$.
Restricting the exact sequence (\ref{sequence}) to $nD$, we obtain
an exact sequence
\begin{equation}
\label{first restriction}
\shTor_1(\O_{nD},\O_C(tA))\lra \shE_1|_{nD}\lra\shE_0|_{nD}\lra
\O_{C\cap nD}\lra 0.
\end{equation}
The map on the left vanishes: on the one hand,
$\shTor_1(\O_{nD},\O_C(tA))$ is supported on
$C\cap D$. On the other hand, $\shE_1|_{nD}$ has no torsion sections
near $C\cap D$, because $D$ is Cartier and $C$ is contained in the
$S_2$-locus of $X$. We conclude that the restriction $\shE_1|_{nD}$
remains an elementary transformation of $\shE_0|_{nD}$.

Next, consider the invertible $\O_C$-module
$\shN_1=\shExt^1(\shL_C(tA),\O_X)$. Dualizing the exact sequence
(\ref{sequence}), we obtain an exact sequence
$$
0\lra \shE_0^\vee\lra\shE_1^\vee\lra \shN_1\lra 0
$$
Its restriction to $C$ gives an exact sequence
\begin{equation}
\label{dual}
0\lra\shN_2\lra\shE_1^\vee|_C\lra \shN_1\lra 0
\end{equation}
for some invertible $\O_C$-module $\shN_2$.
Now choose an ample Cartier divisor
$A'\subset C$ disjoint from
$D\cap C$. Choose $k_0>0$ so that we have
$H^1(C,\O_C(kA'-nD))=0$ and
$H^1(C,\shN_1\otimes\shN_2^\vee((k-1)A'-nD))=0$ for all
$k\geq k_0$.
Fix such an integer $k\geq k_0$, and let $s_1':\shN_1\ra\shN_1(kA')$
be the canonical map given by
$kA'$. The exact sequence (\ref{dual}) yields an exact sequence
$$
H^0(C,\shN_2^\vee\otimes\shN_1(kA'))\lra
H^0(C,\O_{A'}\oplus\O_{C\cap nD})\lra
H^1(\shN_2^\vee\otimes\shN_1((k-1)A'-nD)),
$$
showing that some map $s_2':\shN_2\ra\shN_1(kA')$ has no zeros on
$A'$ and vanishes on $C\cap nD$. The exact sequence
$$
\Hom(\shE_1^\vee|_C,\shN_1(kA'-nD))
\ra\Hom(\shN_2,\shN_1(kA'-nD))\ra\Ext^1(\shN_1,\shN_1(kA'-nD))
$$
tells us that there is an extension of the map
$s_2':\shN_2\ra\shN_1(kA')$ to a mapping
$s_2':\shE_1^\vee|_C\ra\shN_1(kA')$ vanishing on $C\cap nD$.
Setting $s'=s'_1+s_2'$, we obtain a surjective mapping
$s':\shE_1^\vee|_C\ra\shN_1(kA')$. By construction, this
surjection and the old surjection
$\shE_1^\vee|_C\ra\shN_1$ are isomorphic on
$C\cap nD$. Now the exact sequence
$$
0\lra\shE_2\lra\shE_1^\vee\lra\shN_1(kA')\lra 0
$$
defines a locally free $\O_X$-module $\shE_2$.
Set $\shE=\shE_2^\vee$. Taking determinants on $X-A'$, we obtain
$\det(\shE)=\shL$. Using the Riemann--Roch formula (\ref{RR}) on the
$S_2$-ization of $X$, we compute
$$
c_2(\shE)=c_2(\shE_1) + C\cdot\det(\shE_1) + k\deg(A')
+\deg(\shN_1),
$$
which becomes arbitrarily large for $k\gg 0$.
Finally, restriction to $nD$ gives an exact sequence
\begin{equation}
\label{second restriction}
0\lra\shE^\vee|_{nD}\lra\shE^\vee_1|_{nD}\lra \O_{C\cap nD}\lra 0.
\end{equation}
By construction, the elementary transformation
(\ref{first restriction}) and (\ref{second restriction}) are dual to
each other, so $\shE|_{nD}\simeq\shE_0|_{nD}$ is trivial.
Thus $\shE$ is a locally free $\O_X$-module with the desired
properties.
\qed

\begin{remark}
The preceding proof shows that we may increase $c_2(\shF)$ without
changing the determinant by adding the number $\deg(A')>0$.
If $Y$ is reduced and the ground field $k$ is algebraically closed,
we may  choose the Cartier divisor $C\subset X$ reduced and
the ample Cartier divisor $A'\subset C$ of degree one.
In this case, there is a constant $c\geq 0$ such that all numbers
$c_2\geq c$ occur as second Chern number of  locally free
$\O_Y$-modules of rank two with fixed determinant.
\end{remark}

\begin{remark}
If $\shM_1$ and $\shM_2$ and two invertible $\O_Y$-modules, the
locally free $\O_Y$-module
$\shF=\shM_1\oplus\shM_2$ has
$c_2(\shF)=\shM_1\cdot\shM_2$. If the intersection form on $\Pic(Y)$
is nonzero, this immediately implies that $c_2$ can be
arbitrarily positive and negative.
Moreover, if the intersection form has a positive eigenvalue, say
$\shM\cdot\shM>0$, then $\shF=\shM\oplus\shM^\vee$ satisfies
$c_1(\shF)=0$ and
$c_2(\shF)=-\shM\cdot\shM$, so
$c_2$ can be arbitrarily negative with
$c_1$ fixed.

Theorem \ref{existence}
is mainly a result about surfaces whose intersection form
is trivial. However, the intersection form is usually
trivial for singular surfaces, compare \cite{Schroeer 1999}
and \cite{Hartshorne 1977}, Exercise 5.9 on page 232.
\end{remark}

\begin{remark}
\mylabel{complex}
B\u{a}nic\u{a} and Le Potier (\cite{Banica; Le Potier 1987},
Theorem 3.1)
proved
$2c_2(\shF)\geq  c_1^2(\shF)$
for locally free sheaves $\shF$ of rank two on
nonalgebraic smooth compact complex surfaces.
In light of this, we do not expect that,
on badly singular proper algebraic surfaces, the
second Chern numbers for vector bundles of rank two with fixed
determinant  might be arbitrarily negative.
\end{remark}

\section{Existence of global resolutions}

Let $Y$ be a separated surface. According to the
Nagata Embedding Theorem (\cite{Luetkebohmert 1993}, Theorem 3.2),
we may embed $Y$ into a proper surfaces, and Theorem \ref{existence}
tells us that there are many locally free $\O_Y$-modules.
A natural question to ask: Given a coherent $\O_Y$-module
$\shM$, does there exist a locally free $\O_Y$-module
$\shF$ of finite rank,
together with a surjection $\shF\ra\shM$?
This holds if $Y$ admits an ample invertible sheaf, or more
generally an ample family of invertible sheaves
(\cite{SGA 6}, Expos\'e II, Proposition 2.2.3).
Schuster solved the analogous problem for smooth compact complex
surfaces \cite{Schuster 1982}, which are not necessarily
algebraic. We have the following result:

\begin{theorem}
\mylabel{quotient}
Let $Y$ be a separated normal surface.
Then any coherent $\O_Y$-module $\shM$ is the quotient
of a locally free $\O_Y$-module $\shF$ of finite rank.
\end{theorem}

The idea for the proof is  to construct a sheaf of 1-syzygies
$\shS$ for $\shM$, and then obtain $\shF$ as an extension of $\shM$
by
$\shS$. Similarly, we obtain the   sheaf of 1-syzygies as the
extension of an ideal by a suitable locally free sheaf.
The main problem is to choose these sheaves carefully so
that the desired extensions, which always exist locally,
also exist globally.
We start with some preliminary reductions. Throughout,
$Y$ denotes a normal surface.

\begin{proposition}
\mylabel{reflexive}
Let $\shM$ be a coherent $\O_Y$-module. Then there are finitely
many effective Weil divisors $D_1,\ldots,D_j\subset Y$ and a
surjection
$\bigoplus_{i=1}^j\O_Y(-D_i)\ra\shM$ so
that the induced maps $\bigoplus_{i=1}^j\O_Y(-tD_i)\ra\shM$ remain
surjective for all $t\geq 1$.
\end{proposition}

\proof
Since $\shM$ is coherent, it suffices to see that each germ
$s_y\in\shM_y$, $y\in Y$ lies in
the image of a map $\O_Y(-D)\ra\shM$ for some effective
Weil divisor $D\subset Y$ not containing $y$. Choose a representative
$s_V\in\Gamma(V,\shM)$ on some affine open subset
$V\subset Y$. Since $Y$ is separated, the reduced complement $D=Y-U$
is an effective Weil divisor, say with ideal $\shI=\O_Y(-D)$
(\cite{EGA IVd}, Corollary 21.12.7).
According to \cite{EGA I}, Proposition 6.9.17,
we may extend $s_V:\O_V\ra\shM_V$ to a map $s:\shI^m\ra\shM$ for
certain
$m\geq 1$. Now choose some $n\geq 1$ so that $nD$ contains
$\Spec(\O_Y/\shI^m)$. Then $\O_Y(-nD)\subset\shI^m$, and the
composition
$\O_Y(-nD)\ra\shM$ is the desired map.
\qed

\medskip
Next we reduce to the case that our ground field $k$ is infinite:

\begin{proposition}
\mylabel{separable}
Let $k\subset k^s$ be a separable closure, and $Y^s=Y\otimes k^s$
the induced normal surface.
If all coherent $\O_{Y^s}$-module are  quotients of  locally
free $\O_{Y^s}$-modules of finite rank, than the same holds
for all coherent $\O_Y$-modules.
\end{proposition}

\proof
Let $\shM$ be a coherent $\O_Y$-module. Fix a germ $s_y\in\shM_y$,
$y\in Y$. As in the proof of Proposition \ref{reflexive}, there
is an effective Weil divisor $D\subset Y$ not containing $y$
and a map $\O_Y(-D)\ra\shM$ so that the stalk $s_y$ lies in the
image. According to
\cite{EGA IVc}, Theorem 8.5.2, there is a finite separable field
extension
$k\subset k'$ so that on the normal surface
$Y'=Y\otimes k'$, there is a locally free $\O_{Y'}$-module $\shF'$
and a surjection
$\shF'\ra\O_{Y'}(-D')$, where $D'=D\otimes k'$ is the induced
Weil divisor.
Dualizing, we obtain an injection
$f^*(\O_{Y}(D))\subset{\shF'}^\vee$ which is a direct summand near
$f^{-1}(y)$, where $f:Y'\ra Y$ denotes the canonical projection.
This injection corresponds to an injective map
$\O_{Y}(D)\ra f_*({\shF'}^\vee)$, which is a direct summand near $y$.
Setting $\shF=f_*({\shF'}^\vee)^\vee$, we obtain a map
$\shF\ra\O_{Y}(-D)$ which is surjective near $y$.
The image of the composite map $\shF\ra\shM$ contains the germ
$s_y$, and we infer that $\shM$ is the quotient of a locally free
$\O_Y$-module of finite rank.
\qed

\medskip
The following fact from commutative algebra will be useful:

\begin{proposition}
\mylabel{bijective}
Let $\shI $ and $\shN$ be coherent $\O_Y$-modules.
If $\shN$ is reflexive, then the canonical map
$\shHom(\shI^{\vee\vee},\shN)\ra\shHom(\shI,\shN)$ is
bijective.
\end{proposition}

\proof
This is a local problem, so let us assume that $Y=\Spec(A)$ is
affine. Set $I=\Gamma(Y,\shI)$ and $N=\Gamma(Y,\shN)$. According to
\cite{Hartshorne 1994}, Proposition 1.7, there is an exact sequence
$0\ra N\ra L\ra L'$ for certain free
$A$-modules
$L$ and $L'$. This gives an a commutative diagram with exact rows
$$
\begin{CD}
0@>>>\Hom(I^{\vee\vee},N)@>>>\Hom(I^{\vee\vee},L)@>>>
\Hom(I^{\vee\vee},L')\\
@.@VVV@VVV@VVV\\
0@>>>\Hom(I,N)@>>>\Hom(I,L)@>>>
\Hom(I,L').
\end{CD}
$$
By the 5-Lemma, it suffices to treat the case $N=A$.
Then our map in question becomes the canonical map
$(I^{\vee\vee})^\vee\ra I^\vee$. By \cite{Hartshorne 1994},
Corollary 1.6, both
$(I^{\vee\vee})^\vee$ and $I^\vee$ are reflexive.
At each prime ideal $\primid\subset A$ of height one, the module
$I_\primid$ is the direct sum of a free module and a torsion module.
Hence our map is bijective in codimension one, and \cite{Hartshorne
1994}, Theorem 1.12 implies that it is bijective everywhere.
\qed

\medskip
Let us recall the notion of syzygies.
Let $\shM$ be a coherent $\O_Y$-module. A coherent $\O_Y$-module
$\shS$ is called a \emph{sheaf of 1-syzygies} for $\shM$ if for any point
$y\in Y$, there is an integer $m\geq 0$ and an exact sequence
$0\ra\shS_y\ra\O_{Y,y}^{\oplus m}\ra\shM_y\ra 0$.
Our strategy to construct a surjection $\shF\ra\shM$ from a locally
free sheaf $\shF$ of finite rank is first to prove existence of a
suitable  sheaf of 1-syzygies $\shS$, and then to prove the existence of
an extension
$0\ra\shS\ra\shF\ra\shM\ra0$ globalizing the local extensions.
We now have all means to carry out this plan.

\medskip
\emph{Proof of Theorem \ref{quotient}:}
By the Nagata Embedding Theorem (\cite{Luetkebohmert 1993}, Theorem
3.2), we may embed $Y$ into a proper surface, because $Y$ is
separated. Moreover, any coherent sheaf on $Y$ is the
restriction of a coherent sheaf on the compactification
by \cite{EGA I}, Corollary 6.9.11. So let us
assume that
$Y$ is proper. By Proposition \ref{separable}, we may also assume
that the ground field $k$ is separably closed, hence infinite.

Let $\shM$ be a coherent $\O_Y$-module. We seek a locally free
$\O_Y$-module $\shF$ of finite rank, together with a surjection
$\shF\ra\shM$. In light of Proposition \ref{reflexive}, we may assume
that
$\shM=\O_Y(-D)$ for some effective Weil divisor $D\subset Y$.
Moreover, we may replace $D$ by some high multiple and
assume that $K_Y-2D$ is not linearly equivalent to an
effective Weil divisor.

Choose an integer $r\geq 2$ so that we have for each
$y\in\Sing(Y)$ a surjection
$\O_{Y,y}^{\oplus r+2}\ra\O_{Y,y}(-D)$. The exact sequence
\begin{equation}
\label{syzygy}
0\lra\shS_y\lra\O_{Y,y}^{\oplus r+2}\lra\O_{Y,y}(-D)\lra 0
\end{equation}
defines $\O_{Y,y}$-modules $\shS_y$ of rank $r+1$. Our first task is
to extend these local modules to a coherent
$\O_Y$-module
$\shS$
that is locally free on the regular locus. Then $\shS$ would be a
sheaf of 1-syzygies for $\shM=\O_Y(-D)$.

According to \cite{Evans; Griffith 1985}, Theorem 2.14,
for each $y\in\Sing(Y)$ there is an ideal $\shI_y\subset\O_{Y,y}$
and an exact sequence
\begin{equation}
\label{local extension}
0\lra\O_{Y,y}^{\oplus r}\lra\shS_y\lra\shI_y\lra 0.
\end{equation}
Taking  determinants on $\Spec(\O_{Y,y})-\left\{y\right\}$, we see
that
$\shI_y$ is isomorphic to $\O_{Y,y}(D)$ outside $y$.
Hence there are isomorphisms $\shI_y^{\vee\vee}\ra\O_{Y,y}(D)$, and
the cokernel $\shT_y$ for the composite map
$\shI_y\ra\shI_y^{\vee\vee}\ra\O_{Y,y}(D)$ is supported by
$y\in\Spec(\O_{Y,y})$. Then $\shT=\oplus\shT_y$ is a skyscraper sheaf
supported by
$\Sing(Y)$, and we have a global surjection $\O_Y(D)\ra\shT$. Its
kernel
$\shI\subset\O_Y(D)$ is a coherent fractional $\O_Y$-ideal extending
the given local $\O_{Y,y}$-modules $\shI_y$.

Suppose we have a locally free $\O_Y$-module $\shH$ of rank $r$.
Then $\shExt^1(\shI,\shH)$ is a coherent $\O_Y$-module supported
on
$\Sing(Y)$. Choose  trivializations $\shH_y\simeq\O_{Y,y}^{\oplus
r}$. Then the local extensions (\ref{local extension}) correspond to
a section
$e\in H^0(Y,\shExt^1(\shI,\shH))$. The spectral sequence
$H^p(Y,\shExt^q(\shI,\shH))\Rightarrow Ext^{p+q}(\shI,\shH)$
gives an exact sequence
$$
\Ext^1(\shI,\shH)\lra H^0(Y,\shExt^1(\shI,\shH))\lra
H^2(Y,\shHom(\shI,\shH)).
$$
We want to choose the locally free $\O_Y$-module $\shH$ so that the
term on the right vanishes. For then the section $e\in
H^0(Y,\shExt^1(\shI,\shH))$ comes from a global extension
$0\ra\shH\ra\shS\ra\shI\ra0$, and we have constructed a
sheaf of 1-syzygies
$\shS$ for
$\shM$.
The canonical map $\shHom(\O_Y(D),\shH)\ra\shHom(\shI,\shH)$ is
bijective by Proposition
\ref{bijective}. Serre duality gives
$H^2(Y,\shHom(\O_Y(D),\shH))\simeq
H^0(Y,\shH^\vee\otimes\omega_Y(D))^\vee$. So we seek a locally free
$\O_Y$-module
$\shH$ or rank $r$ without nonzero maps $\shH\ra\omega_Y(D)$.

To achieve this, let $f:X\ra Y$ be a resolution of singularities,
and $E\subset X$ the exceptional curve.
According to Remark \ref{stein}, there is an integer
$n\geq 1$ so that the locally free $\O_Y$-modules of rank $r$
correspond to the locally free $\O_X$-modules of rank $r$
with trivial restriction to $nE$.
Choose a Weil divisor $C$ on $X$ with $\O_X(C)\simeq
f^*(\omega_Y(D))^{\vee\vee}$.

Being regular, the proper surface $X$ is projective
(\cite{Hartshorne 1970}, Theorem 4.2).
Choose a very ample invertible $\O_X$-module $\shL$
satisfying both $H^2(X,\shL(K_X-C))=0$ and
$H^1(X,\shL^{\otimes r}(-nE))=0$. Consider the locally free
$\O_X$-module
$\shE_0=\shL^{\oplus r}$ of rank $r$. Then $\Hom(\shE_0,\O_X(C))=0$.
The restriction $\shE_0|_{nE}$ is globally generated. Set $\Gamma=
H^0(nE,\shE_0|_{nE})$. The canonical surjection
$\Gamma\otimes\O_{nE}\ra \shE_0|_{nE}$ yields a morphism
$\varphi:nE\ra\Grass_r(\Gamma)$ into the
Grassmannian of r-dimensional quotients.
Choose a
generic
$r$-dimensional subvector space $\Gamma'\subset \Gamma$.
Here we use that our ground field is infinite.
For each integer $k\geq 0$, let
$G_k\subset\Grass_r(\Gamma)$ be the subscheme of surjections
$\Gamma\ra\Gamma''$
such
that the composition
$\Gamma'\ra\Gamma''$ has rank
$\leq k$. Note that
$G_{r-1}$ is a reduced Cartier divisor, and that
$G_{r-2}$ has codimension four
(see \cite{Arbarello et al. 1985}, Section II.2).

By the dimensional part of Kleiman's Transversality Theorem
(\cite{Kleiman 1974} Theorem 2), the map
$\varphi:nE\ra\Grass_r(\Gamma)$ is disjoint to
$G_{r-2}$ and passes through
$G_{r-1}$ in finitely many points.
The upshot of  this is that the quotient of the canonical map
$ \Gamma'\otimes \O_{nE}\ra \shE_0|_{nE}$ is an invertible sheaf on
some Cartier divisor
$A\subset nE$. Consequently, we have constructed an exact sequence
$$
0 \lra \O_{nE}^{\oplus r}  \lra \shE_0|_{nE}  \lra \O_A  \lra 0,
$$
and $\shL^{\otimes r}_{nE}\simeq\O_{nE}(A)$.
The exact sequence
$$
H^0(X,\shL^{\otimes r})\lra H^0(nE,\shL^{\otimes r}_{nE})\lra
H^1(X,\shL^{\otimes r}(-nE))
$$
shows that $A=H\cap nE$ for some divisor $H\subset X$ with
$\shL^{\otimes r}\simeq\O_X(H)$.
The exact sequence
$$
H^0(H,\shE_0^\vee\otimes\shN)\lra H^0(
A,\shE_0^\vee\otimes\shN_A)\lra H^1(H,\shE_0^\vee\otimes\shN(-A))
$$
tells us that
the surjection $\shE_0|_{nE} \ra\O_A $ extends to a surjection
$\shE_0\ra\shN$ for any sufficiently ample invertible $\O_H$-module
$\shN$. Then the kernel $\shE\subset\shE_0$ is a locally free
$\O_X$-module of rank $r$ that is trivial on $nE$, hence
$\shH=f_*(\shE)$ is locally free of rank $r$ as well.

Suppose there is a nonzero map $\shH\ra\omega_Y(D)$. Then there is
also a nonzero map $\shE\ra\O_X(C)$. To rule this out,
consider the exact sequence
$$
\Hom(\shE_0,\O_X(C))\lra\Hom(\shE,\O_X(C))\lra\Ext^1(\shN,\O_X(C)).
$$
The term on the right sits in an exact sequence
$$
H^1(X,\shHom(\shN,\O_X(C)))\lra \Ext^1(\shN,\O_X(C))\lra
H^0(Y,\shExt^1(\shN,\O_X(C))).
$$
We have $\shHom(\shN,\O_X(C))=0$ and
$\shExt^1(\shN,\O_X(C))\simeq\shN^\vee(C+H)$, which has no global
sections if $\shN$ is sufficiently ample.
In other words, we may choose $\shN$ so that
$\Hom(\shH,\omega_Y(D))=0$.
In turn, the desired global extension
$$
0\lra\shH\lra\shS\lra\shI\lra 0,
$$
exists, and $\shS$ is a   sheaf of 1-syzygies for $\shM=\O_Y(-D)$.

Now we are almost done. The coherent $\O_Y$-module
$\shExt^1(\shM,\shS)$ is a skyscraper sheaf supported by $\Sing(Y)$.
The local extensions from (\ref{syzygy}) define a global section
$e'\in H^0(X,\shExt^1(\shM,\shS))$. As above, we have an
exact sequence
$$
\Ext^1(\shM,\shS)\lra H^0(Y,\shExt^1(\shM,\shS))\lra
H^2(Y,\shHom(\shM,\shS)).
$$
We claim that the term on the right vanishes.
Indeed, the canonical mapping $\shS\ra\shS'$ into the
bidual $\shS'=\shS^{\vee\vee}$ is bijective
outside $\Sing(Y)$,
so the induced map
$H^2(Y,\shHom(\shM,\shS))\ra H^2(Y,\shHom(\shM,\shS'))$ is
bijective as well. Serre duality tells us that
$H^2(Y,\shHom(\shM,\shS'))\simeq \Hom(\shS',\omega_Y(-D))^\vee$.
The exact sequence
$$
0\lra\shH\lra\shS'\lra\O_Y(D)\lra 0
$$
gives an exact sequence
$$
\Hom(\O_Y(D),\omega_Y(-D))\lra
\Hom(\shS',\omega_Y(-D))\lra\Hom(\shH,\omega_Y(-D)).
$$
The term on the right vanishes, because
$\Hom(\shH,\omega_Y(-D))=0$ by construction of $\shH$.
The term on the left vanishes as well, because we have
chosen the Weil divisor $D\subset Y$ from the very beginning so that
$H^0(Y,\omega_Y(-2D))=0$. Summing up, the section $e'\in
H^0(Y,\shExt^1(\shM,\shS))$ comes from a global extension
$$
0\lra \shS\lra\shF\lra\shM\lra 0,
$$
where $\shF$ is the desired locally free $\O_Y$-module of rank
$r+2$.
\qed

\medskip
The following application of Theorem \ref{quotient} will be useful
later:

\begin{corollary}
\mylabel{q-coh}
Let $Y$ be a normal separated surface. Then the following hold:
\renewcommand{\labelenumi}{(\roman{enumi})}
\begin{enumerate}
\item
Any quasi-coherent $\O_{Y}$-module is a
quotient of a locally free $\O_{Y}$-module of the form
$\bigoplus_{\lambda \in L} \shF_{\lambda}$ with each
$\shF_{\lambda}$ locally free of finite rank and $L$
possibly infinite;
\item
If $\shM \rightarrow \shM'$ is a surjection of quasicoherent
$\O_{Y}$-modules with $\shM'$ coherent then there exists a
locally free $\O_{Y}$-module of finite rank $\shF$ and a
homomorphism
$\shF \rightarrow \shM$ such that the composition $\shF\ra\shM'$
remains surjective.
\end{enumerate}
\end{corollary}

\proof
(i) For each coherent submodule
$\shM_\lambda\subset\shM$ of a given quasicoherent $\O_Y$-module
$\shM$, choose a surjection
$\shF_\lambda\ra\shM_\lambda$ from a locally
free sheaf of finite rank. The induced map
$\bigoplus_{\lambda \in L} \shF_{\lambda}\ra\shM$ is surjective,
because $\shM$ is  the union of the $\shM_\lambda$
by \cite{EGA I}, Corollary 6.9.9.

(ii)  Choose a surjection
$\bigoplus_{\lambda \in L} \shF_{\lambda}\ra\shM$
as in (i). Since $\shM'$ is coherent and $Y$ is
quasicompact, there is   a finite subset $L'\subset L$ such that
$\shF=\bigoplus_{\lambda \in L'} \shF_{\lambda}$ surjects
onto $\shM'$.
\qed

\section{Applications to K-groups}

In this section we collect some more or less obvious applications
of Theorem \ref{quotient}. Throughout, $Y$ denotes
a normal separated surface. We start with the following:

\begin{proposition}
\mylabel{finite dimension}
Let $\shM$ be a coherent $\O_Y$-module whose stalks have
finite projective dimension. Then there is an exact
sequence $0\ra\shF_2\ra\shF_1\ra\shF_0\ra\shM\ra 0$ where
$\shF_0,\shF_1, \shF_2$ are locally free $\O_Y$-modules of finite
rank.
\end{proposition}

\proof
By the Auslander--Buchsbaum Formula (\cite{Eisenbud 1995}, Theorem
19.9), the stalks of $\shM$ have projective dimension
$\pd(\shM_y)\leq
\depth(\O_{Y,y})\leq 2$. According to Theorem \ref{quotient}, we have
an exact sequence $\shF_1\ra\shF_0\ra\shM\ra 0$ with
$\shF_0$ and $\shF_1$ locally free of finite rank.
Then the stalks of the kernel $\shF_2$ for $\shF_1\ra\shF_0$ have
projective dimension
$\leq 0$, hence $\shF_2$ is locally free.
\qed

\medskip

We come to $K$-theory.
Let $W_\perf(Y)$ be the additive category of perfect complexes
of quasicoherent $\O_Y$-modules,
and $W_\naive(Y)$ be the additive category of bounded complexes of
locally free $\O_Y$-modules.
We obtain the corresponding derived categories
$D_\perf(Y)$ and $D_\naive(Y)$ by inverting the quasiisomorphisms.

\begin{proposition}
\mylabel{Waldhausen to triangulated}
The natural inclusion
$W_\naive(Y) \subset W_\perf(Y)$ induces an equivalence of
triangulated categories
$D_\naive(Y)\simeq D_\perf(Y)$.
\end{proposition}

\proof
We shall apply Illusie's result  \cite{SGA 6}, Expos\'e II,
Proposition 1.2 (c). In accordance with Illusie's notation, let 
$\underline{S}$ be the Zariski site of $Y$, $S:=Y$,
$\underline{C}$ be the fibered category (over $\underline{S}$) of 
quasicoherent $\O$-modules,
and $\underline{C}_0$ the fibered subcategory of locally free
$\O$-modules of finite rank.
Illusie's result has four assumptions, which take the following
form:

First, the kernel of a surjection between locally free sheaves
of finite rank is locally free of finite rank
(``$\underline{C}_0$ est localement stable par noyau 
d'epimorphisme''), which is clearly true.
Second, any surjection from a quasicoherent sheaf to a locally free
sheaf of finite rank admits locally a section
(``$\underline{C}_0$ est localement relevable''), which is also true.
Third,  if $\shM\ra\shM'$ is a surjection of quasicoherent sheaves on $Y$,
with
$\shM'$ locally free of finite rank, there is a map from a locally
free sheaf of finite rank $\shF\ra\shM$ such that $\shF\ra\shM'$ is
surjective (``$\underline{C}_{0Y}$ est quasi-relevable''), which
holds by Corollary \ref{q-coh}.
Fourth, any coherent sheaf on $Y$ is the quotient of a locally free sheaf
of finite rank (``tout objet de $\underline{C}_{Y}$ qui est de 
$\underline{C}_0$-type fini est de
$\underline{C}_{0Y}$-type fini''), which is Theorem \ref{quotient}.

Hence Illusie's result applies, and tells us that
$D_\naive(Y)\ra D_\perf(Y)$ is an equivalence
of categories.
\qed

\medskip
The category of perfect complexes $W_\perf(Y)$     carries, in a canonical
way, the structure of a \emph{complicial biWaldhausen category},
as explained in \cite{Thomason; Trobaugh 1990}, Section 1. This
additional structure gives rise to the Waldhausen $K$-theory spectrum
$\KK^\perf(Y)$ and  Waldhausen K-groups
$K^\perf_n(Y)=\pi_n(\KK^\perf(Y))$. We also have  the Quillen
$K$-theory spectrum $\KK^Q(Y)$ and  Quillen $K$-groups $K_n^Q(Y)$,
defined via the \emph{exact category} of locally free $\O_Y$-modules of
finite rank. These two approaches give the same result, up to
homotopy:

\begin{theorem}
\mylabel{K-spectra and higher K-groups}
The natural map of  spectra
$\KK^Q(Y)\ra \KK^\perf(Y)$ is a homotopy equivalence.  In
particular, the induced maps of groups
$K^Q_n(Y)\ra K^\perf_n(Y)$
are bijective for all $n \geq 0$.
\end{theorem}

\proof
Let $\KK^\naive(Y)$ be the Waldhausen $K$-theory spectrum for the
complicial biWaldhausen category $W_\naive(Y)$.
By \cite{Thomason; Trobaugh 1990},
Proposition $3.10$, there is  a canonical homotopy
equivalence
$\KK^Q(Y)\ra\KK^\naive(Y)$. This holds for arbitrary schemes.
By \cite{Thomason; Trobaugh 1990},
Theorem 1.9.8, together with Proposition
\ref{Waldhausen to triangulated}, the   map
$\KK^\naive(Y)\ra \KK^\perf(Y)$ is a homotopy equivalence as well.
\qed

\medskip
The group of connected components for $\KK^Q(X)$ is nothing
but the Grothendieck group for the exact category of
vector bundles. Similarly, the group of connected components
for $\KK^\perf(Y)$ equals the Grothendieck group for
the triangulated category of perfect complexes, which is defined in
\cite{SGA 6}, Expos\'e IV. Therefore, we have the following:

\begin{corollary}
\mylabel{K-group}
The  group $K^{\mathrm{perf}}_{0}(Y)$ is
isomorphic to the free abelian group generated by
isomorphism classes $[\shF]$ of locally free $\O_Y$-modules of
finite rank, modulo the relations $[\shF']+[\shF'']=[\shF]$ for any
exact sequence
$0\ra\shF'\ra\shF\ra\shF''\ra 0$.
\end{corollary}

\begin{remark}
Note that the proof of Corollary \ref{q-coh}
    works on every noetherian scheme $X$ for which any coherent
$\O_{X}$-module is a quotient of a locally free
$\O_{X}$-module of finite rank.
Then, the same argument as in the proof of
Proposition \ref{Waldhausen to triangulated} shows
that for such schemes the Waldhausen $K$-theory of
perfect complexes coincides with the Quillen $K$-theory of vector bundles.
\end{remark}

\section{Counterexamples}

Here we show that, on nonseparated schemes,
there are usually many coherent sheaves that
are not quotients of locally free sheaves.
Let $Z=\Spec(A)$ be an affine normal local noetherian scheme of
dimension at least two, and $U\subset Z$ the complement of the
closed point. Let $Y_1,Y_2$ be two copies of $Z$.
The cocartesian diagram
$$
\begin{CD}
U@>>> Y_1\\
@VVV @VVV\\
Y_2@>>> Y
\end{CD}
$$
defines a nonseparated normal noetherian scheme $Y$,
endowed with an affine covering $Y=Y_1\cup Y_2$
with $Y_1\cap Y_2=U$.

\begin{proposition}
\mylabel{vector bundles}
Every locally free $\O_Y$-module of finite rank is trivial.
\end{proposition}

\proof
Let $\foU$ be the open affine covering $Y=Y_1\cup Y_2$.
We have $Z^1(\foU,\GL_{r})=\Gamma(U,\GL_{r})=\GL_r(A)$ and
$C^0(\foU,\GL_{r})=\prod_{i=1}^2\GL_r(A)$.
This clearly implies that $H^1(\foU,\GL_{r})=0$. Since any locally
free
$\O_Y$-module of finite rank becomes trivial on the local schemes
$Y_i$, we also have
$H^1(Y,\GL_{r})=0$.
\qed

\begin{proposition}
\mylabel{vector bundles II}
There are coherent $\O_Y$-modules that are not
quotients of locally free $\O_Y$-modules of finite rank.
\end{proposition}

\proof
Let $\shM$ be a coherent $\O_Y$-module.
Suppose there is an exact sequence
$$
\O_Y^{\oplus r}\lra \O_Y^{\oplus s}\lra \shM\lra 0.
$$
The map on the left is given by an $s\times r$-matrix with entries
in $\Gamma(Y,\O_Y)=A$. This implies $\shM|_{Y_1}\simeq\shM|_{Y_2}$,
where we use the canonical identifications $Y_1=Z=Y_2$.
Now choose a coherent $\O_Y$-module $\shM$ so that
$\shM|_{Y_1}\not\simeq\shM|_{Y_2}$, for example with $\shM_U=0$. Then
either there is no surjection
$\O_Y^{\oplus s}\ra \shM$, or there is such a surjection and
its kernel $\shM'$ does not admit a surjection $\O_Y^{\oplus r}\ra
\shM'$.
The statement now follows from Proposition \ref{vector bundles}.
\qed


\end{document}